\input amssym.def
\input amssym.tex


\catcode`\@=11

\hsize=125 mm   \vsize =187mm
\hoffset=4mm    \voffset=10mm
\pretolerance=500 \tolerance=1000 \brokenpenalty=5000

\catcode`\;=\active
\def;{\relax\ifhmode\ifdim\lastskip>\z@
\unskip\fi\kern.2em\fi\string;}

\catcode`\:=\active
\def:{\relax\ifhmode\ifdim\lastskip>\z@\unskip\fi
\penalty\@M\ \fi\string:}

\catcode`\!=\active
\def!{\relax\ifhmode\ifdim\lastskip>\z@
\unskip\fi\kern.2em\fi\string!}

\catcode`\?=\active
\def?{\relax\ifhmode\ifdim\lastskip>\z@
\unskip\fi\kern.2em\fi\string?}

\def\^#1{\if#1i{\accent"5E\i}\else{\accent"5E #1}\fi}
\def\"#1{\if#1i{\accent"7F\i}\else{\accent"7F #1}\fi}


\catcode`\@=12

\newif\ifpagetitre      \pagetitretrue
\newtoks\hautpagetitre  \hautpagetitre={\hfil}
\newtoks\baspagetitre   \baspagetitre={\hfil}

\newtoks\auteurcourant  \auteurcourant={\hfil}
\newtoks\titrecourant   \titrecourant={\hfil}

\newtoks\hautpagegauche \newtoks\hautpagedroite
\hautpagegauche={\hfil\the\auteurcourant\hfil}
\hautpagedroite={\hfil\the\titrecourant\hfil}

\newtoks\baspagegauche  \baspagegauche={\hfil\tenrm\folio\hfil}
\newtoks\baspagedroite  \baspagedroite={\hfil\tenrm\folio\hfil}

\headline={\ifpagetitre\the\hautpagetitre
\else\ifodd\pageno\the\hautpagedroite
\else\the\hautpagegauche\fi\fi}
\footline={\ifpagetitre\the\baspagetitre
\global\pagetitrefalse
\else\ifodd\pageno\the\baspagedroite
\else\the\baspagegauche\fi\fi}

\hautpagetitre={\hfill\tenbf preliminary version: Not for diffusion\hfill}
\hautpagetitre={\hfill}
\hautpagegauche={\tenrm\folio\hfill\the\auteurcourant}
\hautpagedroite={\tenrm\the\titrecourant\hfill\folio}
\baspagegauche={\hfil} \baspagedroite={\hfil}
\auteurcourant{Marmi, Moussa, Yoccoz}
\titrecourant{On the cohomological equation for interval exchange maps}
\def\mois{\ifcase\month\or January\or February\or March\or April\or
May\or June\or July\or August\or September\or October\or November\or
December\fi}
\def\Date{\rightline{\mois\ /\ \the\day\ /\/ \the\year}}
\hfuzz=0.3pt
\font\tit=cmb10 scaled \magstep1

\def\H{\Bbb H}
\def\R{\Bbb R}

\def\Z{\Bbb Z}

\def\hH+{\hat{\H}^{+}}
\def\hHZ+{\widehat{\H^{+}/\Z}}

\def\PGL2Z{\hbox{PGL}\, (2,\Z)}
\def\GL2Z{\hbox{GL}\, (2,\Z)}
\def\SL2Z{\hbox{SL}\, (2,\Z)}

\def\cA{{\cal A}}

\def\cR{{\cal R}}

\def\cV{{\cal V}}

\def\Proc#1#2\par{\medbreak \noindent {\bf #1\enspace }{\sl #2}%
\par\ifdim \lastskip <\medskipamount \removelastskip \penalty 55\medskip \fi}%
\def\Def#1#2\par{\medbreak \noindent {\bf #1\enspace }{#2}%
\par\ifdim \lastskip <\medskipamount \removelastskip \penalty 55\medskip
\fi}

\def\hfl#1#2{\smash{\mathop{\hbox to 12mm{\rightarrowfill}}
\limits^{\scriptstyle#1}_{\scriptstyle#2}}}

\def\build#1_#2^#3{\mathrel{\mathop{\kern 0pt#1}\limits_{#2}^{#3}}}



\font\pic=cmr7

\rightline{March 17, 2003}
\medskip
\centerline{\tit On the cohomological equation 
for interval exchange maps}
\bigskip
\centerline{\tit Sur l'\'equation cohomologique 
pour les \'echanges d'intervalles}
\bigskip
\centerline{S. Marmi\footnote{$^1$}{Dipartimento di Matematica e
Informatica,
Universit\`a di Udine, Via delle Scienze 206, Loc. Rizzi, 33100 Udine,
Italy
and Scuola Normale Superiore, Piazza dei Cavalieri 7, Pisa, Italy},
P. Moussa\footnote{$^2$}{Service de Physique Th\'eorique, CEA/Saclay,
91191 Gif-Sur-Yvette, France}
and J.-C. Yoccoz\footnote{$^3$}{Coll\`ege de France, 3, Rue d'Ulm,
75005 Paris and
Universit\'e de Paris-Sud,
Math\'ematiques. B\^at. 425, 91405-Orsay, France}}
\vskip 2. truecm
\centerline{\bf  Abstract}

\vskip .5 truecm
\noindent
We exhibit an explicit full measure class of minimal interval exchange maps 
$T$ for which the cohomological equation $\Psi -\Psi\circ T=\Phi$
has a bounded solution $\Psi$ provided that the datum $\Phi$
belongs to a finite codimension subspace of the space of functions 
having on each interval a derivative of bounded variation.

\noindent
The class of interval exchange maps is characterized in terms
of a diophantine condition of ``Roth type'' imposed to an 
acceleration of the Rauzy--Veech--Zorich continued fraction 
expansion associated to $T$.

\vskip 1. truecm
\centerline{\bf  R\'esum\'e}

\vskip .5 truecm
\noindent
On pr\'esente une classe explicite d'\'echanges 
d'intervalles $T$, de mesure pleine, pour laquelle 
l'\'equation cohomologique $\Psi -\Psi\circ T=\Phi$
admet une solution born\'ee $\Psi$, \`a condition que la 
donn\'ee $\Phi$ appartienne \`a un sous--espace de codimension
finie de l'espace des fonctions dont la d\'eriv\'ee sur 
chaque intervalle est de variation born\'ee.

\noindent
Cette classe est d\'efinie par une condition diophantienne 
``de type Roth'' exprim\'e dans une variante du 
d\'eveloppement en fraction continue de Rauzy--Veech--Zorich
associ\'e \`a $T$.

\vskip 1. truecm
\centerline{CONTENTS}
\smallskip{\pic

\item{0.} French abridged version
\item{1.} Interval exchange maps and the cohomological equation. 
Main Theorem
\item{2.} Rauzy--Veech--Zorich continued fraction algorithm 
and its acceleration
\item{3.} Special Birkhoff sums
\item{4.} The Diophantine condition
\item{5.} Sketch of the proof of the theorem
}

\vskip 1. truecm
{\it Acknowledgements.} This work began in 1998. This research has been
supported by the CNR, CNRS, INDAM--GNFM and MURST.

\vfill\eject
\noindent
{\bf 0. French abridged version}

\vskip .5 truecm
Soit $\cA$ un alphabet constitu\'e de $d\ge 2$ lettres.
Soit $(\pi_0, \pi_1)$ un couple de {\it bijections } de 
$\cA$ sur $\{1, \ldots ,d\}$, qu'on supposera 
toujours {\it admissible}: pour tout $1\le k<d$,
on a 
$\pi_0^{-1}(\{1,\ldots ,k\})\not= \pi_1^{-1}(\{1,\ldots ,k\})\,$.
Si on se donne aussi des longueurs $(\lambda_\alpha)_{\alpha\in\cA}$
(avec $\lambda_\alpha>0$), on d\'efinit un \'echange d'intervalles 
$T$ [V1] (avec discontinuit\'es marqu\'ees) par la formule 
suivante: pour $\alpha\in\cA$, $0\le x<\lambda_\alpha$
$$
T\left(x+\sum_{\pi_0 (\beta )<\pi_0 (\alpha )}\lambda_\beta \right)
=x+\sum_{\pi_1 (\beta )<\pi_1 (\alpha )}\lambda_\beta\; .
$$
On posera 
$I_\alpha =\left[\sum_{\pi_0 (\beta )<\pi_0 (\alpha )}\lambda_\beta ,
\sum_{\pi_0 (\beta )\le\pi_0 (\alpha )}\lambda_\beta \right)$.
Comme $T$ est une translation sur chaque $I_\alpha$, $T$ pr\'eserve 
l'orientation et la mesure de Lebesgue. On peut penser $\cA$ 
comme l'ensemble des orbites des points de discontinuit\'e de $T$.

\noindent
\vskip .3 truecm
Notons $\Delta$ le simplexe standard de dimension $d-1$ dans $\R^d$.
Alors le vecteur normalis\'e des longueurs $\left(\lambda_\alpha
\left(\sum_{\beta\in\cA}\lambda_\beta\right)^{-1}\right)_{\alpha\in\cA}$
appartient \`a  $\Delta$: on notera $\Delta (\pi_0 ,\pi_1)$ 
le simplexe des longueurs normalis\'ees associ\'e \`a un couple admissible
$(\pi_0 ,\pi_1)$.

On notera $\hbox{BV}\, (\sqcup I_\alpha)$ 
l'espace des fonctions $\varphi$ sur $[0,1)$ dont la restriction \`a
chaque $I_\alpha$ est de variation born\'ee, et 
$\hbox{BV}_*\, (\sqcup I_\alpha)$ l'hyperplan form\'e des 
fonctions de moyenne nulle.

\vskip .3 truecm\noindent
{\bf Th\'eor\`eme.}{\it Soit $(\pi_0 ,\pi_1)$ un couple admissible. Il existe
une partie $D(\pi_0 ,\pi_1)\subset \Delta (\pi_0 ,\pi_1)$ de {\rm mesure totale}
telle que, pour tout $T\in D(\pi_0 ,\pi_1)$, tout $\varphi
\in \hbox{\rm BV}_*\, (\sqcup I_\alpha)$ on puisse trouver 
\item{(1)} une fonction 
$\Phi$, Lipschitzienne sur chaque $I_\alpha$, dont la d\'eriv\'ee sur chaque 
$I_\alpha$ est $\varphi$ et de moyenne nulle sur $\sqcup I_\alpha$;
\item{(2)} une fonction born\'ee  $\Psi$, de moyenne nulle sur 
$\sqcup I_\alpha$, 

\noindent
v\'erifiant l'\'equation cohomologique
$\Psi -\Psi\circ T=\Phi$.
}

\vskip .3 truecm\noindent
{\bf Remarques. } $(1)$ La partie $D(\pi_0 ,\pi_1)$ of $\Delta 
(\pi_0 ,\pi_1)$ est explicitement d\'efinie ci--dessous.
\item{\phantom{1.3 Remarks. } $(2)$}
Keane [Ke1] a montr\'e que si aucune orbite de $T$ ne passe plus 
d'une fois par une discontinuit\'e de $T$, alors $T$ est minimal.
Les \'echanges d'intervalles dans $D(\pi_0
,\pi_1)$  satisfont l'hypoth\`ese de Keane. 
Masur [Ma] et Veech [V2] ont montr\'e que presque
tout \'echange  d'intervalles est 
uniquement ergodique (pour $d\ge 4$, il existe des \'echanges d'intervalles qui ne sont pas 
uniquement ergodiques [KN, Ke2]). Les \'echanges d'intervalles dans 
$D(\pi_0 ,\pi_1)$ sont uniquement ergodiques.
\item{\phantom{1.3 Remarks. } $(3)$} Notre r\'esultat est clairement 
reli\'e (via une suspension singuli\`ere) au th\'eor\`eme de Forni
[Fo] sur l'\'equation cohomologique pour les champs de vecteurs 
pr\'eservant les aires sur les surfaces. Notre m\'ethode est diff\'erente
de celle de Forni: il construit des instruments d'analyse de Fourier sur les 
surfaces plates (singuli\`eres) et utilise le th\'eor\`eme de Fatou
sur les valeurs au bord des fonctions holomorphes born\'ees. Nous 
mettons \`a profit un algorithme de fraction
continue d\^u \`a Rauzy--Veech--Zorich pour
obtenir une condition diophantienne explicite. La
perte de diff\'erentiabilit\'e est par ailleurs
plus faible que  dans son th\'eor\`eme.

\vskip .3 truecm\noindent
{\bf Esquisse de preuve.}

\noindent
Soit $T=T(0)$ un \'echange d'intervalles v\'erifiant 
l'hypoth\`ese de Keane. Une version acc\'el\'er\'ee de 
l'algorithme de Rauzy--Veech--Zorich ([Ra], [V2], [Z1], 
[Z2]) fournit une suite d'\'echanges d'intervalles 
$(T(n))_{n\ge 0}$ avec des longueurs associ\'ees $(\lambda_\alpha
(n))_{\alpha\in\cA\, , \,n\ge 0}$ qui satisfont aux propri\'et\'es 
suivantes:
\item{(i)} pour $m<n$, $T(n)$ est l'application de premier 
retour de $T(m)$ sur $I(n)=[0,\sum_{\alpha\in\cA}\lambda_\alpha
(n))$;
\item{(ii)} pour $m<n$, on a $\lambda (m) =Q(m,n)\lambda (n)$,
avec une matrice $Q(m,n)\in\hbox{SL}\, (d,\Z )$ \`a coefficients $\ge 0$;
\item{(iii)} pour $m<n$,$\alpha ,\beta\in\cA$, le temps pass\'e par 
$I_\beta (n)$ dans $I_\alpha (m)$ avant retour dans $I(n)$ est 
$Q_{\alpha\beta}(m,n)$, et le temps de retour est $Q_\beta (m,n)=
\sum_{\alpha\in\cA}Q_{\alpha\beta}(m,n)$;
\item{(iv)} on a $Q_{\alpha\beta}(m,n)>0$ pour $\alpha ,\beta\in\cA$, 
$n\ge m+m(d)$.

\noindent
Soit $\varphi$ une fonction sur $\sqcup_{\alpha\in\cA}I_\alpha (m)$; 
pour $n>m$, $\beta\in\cA$, $x\in I_\beta (n)$, on d\'efinit
$$
S(m,n)\varphi (x) = \sum_{0\le l<Q_\beta (m,n)}
\varphi((T(m))^l(x))\; .
$$
Notons $\Gamma^{(m)}$ l'espace de dimension $d$ des fonctions constantes
sur chaque $I_\alpha (m)$; alors $S(m,n)$ envoie $\Gamma^{(m)}$ sur 
$\Gamma^{(n)}$, sa matrice par rapport aux bases canoniques \'etant
${}^tQ(m,n)$. On posera
$$
\Gamma^{(m)}_s=\{\chi\in \Gamma^{(m)}\, , \,
\limsup_{n\to\infty} {\log\Vert S(m,n)\chi\Vert\over\log
\Vert Q(m,n)\Vert} <0\}\; .
$$
On dira que $T$ est de type Roth si:
\item{(a)} pour tout $\varepsilon >0$, on a $\Vert Q(n,n+1)\Vert \le 
\Vert Q(0,n)\Vert^\varepsilon$ pour tout $n$ assez grand;
\item{(b)} Il existe $\theta >0$ tel qu'on ait 
$\Vert S(0,n)\vert_{\Gamma^{(0)}_*}\Vert\le \Vert Q(0,n)
\Vert^{1-\theta}$ pour tout $n$ assez grand, $\Gamma^{(0)}_*$
d\'esignant l'hyperplan de $\Gamma^{(0)}$ form\'e des fonctions 
de moyenne nulle;
\item{(c)} L'op\'erateur $S(m,n)_\flat\, :\, \Gamma^{(m)}/\Gamma^{(m)}_s
\rightarrow \Gamma^{(n)}/\Gamma^{(n)}_s$ induit par $S(m,n)$ v\'erifie, pour
tout $\varepsilon >0$, $\Vert (S(m,n)_\flat )^{-1}\Vert\le
\Vert Q(0,n)\Vert^\varepsilon$ si $n$ est assez grand.

\noindent
Sous les conditions (a) et (b), on montre que, pour une fonction 
\`a variation born\'ee $\varphi$ sur $\sqcup_{\alpha\in\cA}I_\alpha (0)$, on
a
$$
\Vert S(0,n)(\varphi )\Vert_{L^\infty}\le \Vert Q(0,n)\Vert^{1-\theta'}
\Vert\varphi\Vert_{BV}\; ,
$$
avec $\theta'=\theta'(\theta )>0$, pour $n$ assez grand. En utilisant aussi
(c), on arrive \`a construire une fonction lipschitzienne $\Phi$ dont la 
d\'eriv\'ee sur chaque $I_\alpha (0)$ est $\varphi$, qui v\'erifie
$$
\sum_{n\ge 0}\Vert Q(n,n+1)\Vert\Vert S(0,n)(\Phi )\Vert_{L^\infty}
<+\infty\; .
$$
On en d\'eduit que les sommes de Birkhoff de $\Phi$ sont born\'ees et on 
conclut gr\^ace \`a un th\'eor\`eme de Gottschalk--Hedlund.

\vskip .5 truecm
\noindent
{\bf 1. Interval exchange maps and the cohomological equation. Main Theorem}

\vskip .5 truecm
{\bf 1.1} Let $\cal A$ denote an alphabet with $d\ge 2$ elements. 
Consider a pair $(\pi_0, \pi_1)$ of {\it bijections } of 
$\cA$ on $\{1, \ldots ,d\}$:
we will always assume that the pair is {\it admissible}: for all
$1\le k<d$ one has 
$\pi_0^{-1}(\{1,\ldots ,k\})\not= \pi_1^{-1}(\{1,\ldots ,k\})\,$.
One can associate the permutation $\pi=\pi_1\circ\pi_0^{-1}$
of $\{1,\ldots ,d\}$ to the pair $(\pi_0, \pi_1)$.
Given an admissible pair $(\pi_0, \pi_1)$
and a vector of lengths of intervals $(\lambda_\alpha)_{\alpha\in\cA}$,
an {\it interval exchange map} (i.e.m.\ ) $T$ [V1]({\it with marked 
discontinuities}) is defined through the formula: for $\alpha\in\cA$, 
$0\le x<\lambda_\alpha$
$$
T\left(x+\sum_{\pi_0 (\beta )<\pi_0 (\alpha )}\lambda_\beta \right)
=x+\sum_{\pi_1 (\beta )<\pi_1 (\alpha )}\lambda_\beta\; .
\eqno(1)
$$
We write 
$I_\alpha =\left[\sum_{\pi_0 (\beta )<\pi_0 (\alpha )}\lambda_\beta ,
\sum_{\pi_0 (\beta )\le\pi_0 (\alpha )}\lambda_\beta \right)$.
We can think of $\cA$ as the set of the orbits 
of the discontinuities, taking right limits 
(each interval is associated to the discontinuity at its left endpoint). 
Each i.e.m.\ is  piecewise a translation, orientation--preserving and 
preserves Lebesgue measure. 

\smallskip
{\bf 1.2}
Let $\Delta$ denote the standard $(d-1)$--dimensional simplex
in $\R^d$. The normalized vector $\left(\lambda_\alpha
\left(\sum_{\beta\in\cA}\lambda_\beta\right)^{-1}\right)_{\alpha\in\cA}$
belongs to $\Delta$: we denote $\Delta (\pi_0 ,\pi_1)$ 
the symplex of normalized lengths corresponding to the admissible 
pair $(\pi_0 ,\pi_1)$.
We will denote $\hbox{BV}\, (\sqcup I_\alpha)$ 
(resp.\ $\hbox{BV}_*\, (\sqcup I_\alpha)$) 
the space of functions $\varphi$ whose restriction to each of the intervals 
$I_\alpha$ is a function of bounded variation
(resp.\ the hyperplane of 
$\hbox{BV}\, (\sqcup I_\alpha)$ made of functions
whose integral on the disjoint union $\sqcup I_\alpha$ vanishes).

Our main result can be stated as follows: 

\vskip .3 truecm\noindent
{\bf Theorem.}{\it Let $(\pi_0 ,\pi_1)$ be admissible. There exists a 
subset $D(\pi_0 ,\pi_1)\subset \Delta (\pi_0 ,\pi_1)$ of {\rm full measure}
such that for all $T\in D(\pi_0 ,\pi_1)$ and for all function $\varphi
\in \hbox{\rm BV}_*\, (\sqcup I_\alpha)$ one can find
\item{(1)} a function 
$\Phi$, Lipschitz on each $I_\alpha$, with $D\Phi=\varphi$ 
on each $I_\alpha$ and total mean value $0$ on $\sqcup I_\alpha$;
\item{(2)} a bounded function $\Psi$, 
with total mean value $0$ on $\sqcup I_\alpha$, 

\noindent
which satisfy the 
cohomological equation $\Psi -\Psi\circ T=\Phi$.
}

\vskip .3 truecm\noindent
{\bf 1.3 Remarks. } $(1)$ The subset $D(\pi_0 ,\pi_1)$ of $\Delta 
(\pi_0 ,\pi_1)$ will be explicitely defined below.
\item{\phantom{1.3 Remarks. } $(2)$}
Keane [Ke1] proved that if there is no orbit segment starting and ending
with discontinuities, $T$ is minimal. All the i.e.m.\ considered 
here are assumed to satisfy this hypothesis, which in particular holds
when the lengths are rationally independent. 
Masur [Ma] and Veech [V2] proved that almost all i.e.m.\
are uniquely ergodic (for $d\ge 4$, there exist non 
uniquely ergodic i.e.m.\ [KN, Ke2]).
The diophantine condition $D(\pi_0 ,\pi_1)$ is easily seen to imply
unique ergodicity.
\item{\phantom{1.3 Remarks. } $(3)$} Obviously, our result is 
closely connected (via singular suspension) to Forni's theorem
[Fo] on the cocycle equation for area--preserving vector fields on surfaces.
Our method is different from Forni's: he constructs some 
Fourier analysis on flat surfaces and gets his results 
through Fatou's theorem on boundary values of bounded 
holomorphic functions. Our method provides an explicit diophantine 
condition and a smaller loss of differentiability.
\item{\phantom{1.3 Remarks. } $(4)$}  We can prove
similar results, with the same loss of
differentiability, when we start with a more
regular data $\Phi$.
\vskip .5 truecm
\noindent
{\bf 2.
Rauzy--Veech--Zorich continued fraction algorithm 
and its acceleration.}

\vskip .5 truecm
In order to describe the set $D(\pi_0,\pi_1)$ we will make use of 
the generalization of continued fractions to i.e.m.\ 's due to the 
work of Rauzy [Ra], Veech [V2] and Zorich [Z1,Z2].

\vskip .3 truecm\noindent
{\bf 2.1}
Let $(\pi_0, \pi_1)$ be an admissible pair. We define two 
new admissible 
pairs $\cR_0 (\pi_0, \pi_1)$ and $\cR_1 (\pi_0, \pi_1)$ as follows: 
let $\alpha_0, \alpha_1$ be the (distinct) elements of $\cA$ such that 
$\pi_0(\alpha_0)=\pi_1(\alpha_1)=d$; one has
$$
\eqalign{
\cR_0 (\pi_0, \pi_1) &= (\pi_0, \hat{\pi}_1)\; , \cr
\cR_1 (\pi_0, \pi_1) &= (\hat{\pi}_0, \pi_1)\; , \cr
}\eqno(2)
$$
where
$$
\eqalign{
\hat{\pi}_1 (\alpha) &=\cases{
\pi_1(\alpha ) & if $\pi_1(\alpha )\le \pi_1(\alpha_0)$,\cr
\pi_1(\alpha )+1 & if $\pi_1(\alpha_0 )<\pi_1(\alpha)<d$,\cr
\pi_1(\alpha_0 )+1 & if $\alpha=\alpha_1$, ($\pi_1(\alpha_1)=d$);\cr
}\cr
\hat{\pi}_0 (\alpha) &=\cases{
\pi_0(\alpha ) & if $\pi_0(\alpha )\le \pi_0(\alpha_1)$,\cr
\pi_0(\alpha )+1 & if $\pi_0(\alpha_1 )<\pi_0(\alpha)<d$,\cr
\pi_0(\alpha_1 )+1 & if $\alpha=\alpha_0$, ($\pi_0(\alpha_0)=d$).\cr
}\cr
}\eqno(3)
$$

\noindent
The {\it extended Rauzy class} of $(\pi_0, \pi_1)$ is the 
set of admissible 
pairs obtained by saturation of $(\pi_0, \pi_1)$ under the action of
$\cR_0$ and $\cR_1$. 
The {\it extended Rauzy diagram} has for vertices the elements
of the extended Rauzy class, each vertex $(\pi_0,\pi_1)$
being the origin of two arrows joining $(\pi_0,\pi_1)$
to $\cR_0 (\pi_0, \pi_1)$, $\cR_1 (\pi_0, \pi_1)$. The
{\it name} of an arrow joining $(\pi_0,\pi_1)$
to $\cR_\varepsilon (\pi_0, \pi_1)$ (with 
$\varepsilon\in\{0,1\}$) is the element $\alpha_\varepsilon\in
\cA$ such that $\pi_\varepsilon (\alpha_\varepsilon )=1$.

\vskip .3 truecm\noindent
{\bf 2.2}
Let $T$ be an i.e.m.\ with marked discontinuities,
given by data $(\pi_0, \pi_1)$, $(\lambda_\alpha )_{\alpha
\in\cA}$.
For $\varepsilon\in\{0,1\}$, define $\alpha_\varepsilon\in\cA$
by $\pi_\varepsilon (\alpha_\varepsilon )=d$  as above.

We say that $T$ is of {\it type $\varepsilon$} 
if one has $\lambda_{\alpha_\varepsilon}\ge\lambda_{
\alpha_{1-\varepsilon}}$;
we then define a new i.e.m.\ $\cV (T)$ by the following data: 
the admissible pair $\cR_\varepsilon(\pi_0, \pi_1)$ and the lengths
$(\hat{\lambda}_\alpha)_{\alpha\in\cA}$ given by
$$
\cases{
\hat{\lambda}_\alpha = \lambda_\alpha & if $\alpha\not=\alpha_\varepsilon$,\cr
\hat{\lambda}_{\alpha_\varepsilon} = \lambda_{\alpha_\varepsilon}-
\lambda_{\alpha_{1-\varepsilon}} & otherwise. \cr
}\eqno(4)
$$
The i.e.m.\ $\cV (T)$ is the first return map of  $T$ on 
$\left[0,\sum_{\alpha\not=\alpha_{1-\varepsilon}}\lambda_\alpha\right)$.
We also associate to $T$ the arrow in the extended Rauzy diagram 
joining $(\pi_0,\pi_1)$ to $\cR_\varepsilon(\pi_0, \pi_1)$.
Iterating this process, we obtain a sequence of i.e.m\ 
$(\cV^k (T))_{k\ge 0}$ and an infinite path in the extended 
Rauzy diagram starting from $(\pi_0,\pi_1)$.

\vskip.3 truecm\noindent
{\bf 2.3}
The next step is to group together several iterations of
$\cV$ to obtain the Zorich and the accelerated Zorich algorithms.

Starting from $T=T(0)$, we define a sequence $T(n)=
\cV^{k(n)}(T)$ by the following property: for $n\ge 0$, 
$k(n+1)$ is the largest integer $k>k(n)$ such that not all names
in $\cA$ are taken by arrows associated to iterations of
$\cV$ from $T(n)$ to $\cV^k(T)$. This definition is based on 
the following elementary lemma:

\vskip .3 truecm\noindent
{\bf Lemma}{\it Assuming that no orbit segment of $T$ starts and
ends with a discontinuity, every name is taken infinitely many times
in the infinite path associated to $T$.}

\vskip .3 truecm\noindent
{\bf Remark:} We have defined above the accelerated Zorich algorithm,
which is most convenient for us because of the lemma below.
For the Zorich algorithm itself, we have $\tilde T(n)=
\cV^{\tilde k(n)}(T)$, $\tilde k(n+1)$ being the largest integer
$\tilde k>\tilde k(n)$ such that all arrows associated to the iterations 
of $\cV$ from $\tilde T(n)$ to $\cV^k(T)$ have the same name.

\vskip.3 truecm\noindent
{\bf 2.4} Let $T=T(0)$ be an i.e.m.\ satisfying the hypotheses of the lemma
above. Let $(T(n))_{n\ge 0}$ be the sequence of i.e.m.\ obtained 
by the accelerated Zorich algorithm, with associated lengths
$(\lambda_\alpha (n))_{\alpha\in\cA}$.

Iterating formula (4) gives a matrix $Z(n)\in\hbox{SL}\,
(d,\Z )$ with non negative entries such that 
$$
\lambda (n-1)=Z(n)\lambda (n)\; .\eqno(5)
$$
We will write, for $m<n$
$$
\eqalignno{
Q(m,n) &=Z(m+1)\cdots Z(n)\; , &(6)\cr
\lambda (m) &= Q(m,n)\lambda (n)\; . &(7)\cr
}
$$
For $m<n$, $T(n)$ is the first return map of $T(m)$
on $I(n)=\left[0,\sum_{\alpha\in\cA}\lambda_\alpha (n)\right)$;
the return time of $I_\beta (n)$ in $I(n)$ is 
$\sum_\alpha Q_{\alpha\beta}(m,n)$ 
and the time spent in $I_\alpha (m)$ is $Q_{\alpha\beta}(m,n)$.
From the definition of the accelerated Zorich algorithm, one easily gets 
the following:

\vskip .3 truecm\noindent
{\bf Lemma}{\it There exists an integer $m(d)>0$ such that for 
$n\ge m+m(d)$, one has $Q_{\alpha\beta}(m,n)>0$
for all $\alpha ,\beta\in\cA$.}

\vskip .5 truecm
\noindent
{\bf 3. Special Birkhoff Sums}

\vskip .3 truecm\noindent
{\bf 3.1}
Let $T$ as above, $m\le n$, and 
$\varphi$ be a function defined on the disjoint union 
$\sqcup_\cA I_\alpha (m)$. We define a new function 
$S(m,n)\varphi$ on $\sqcup_\cA I_\beta (n)$
as follows: for $x\in I_\beta (n)$
$$
S (m,n)\varphi (x) = \sum_{0\le l<Q_\beta (m,n)}
\varphi (T(m)^l(x)) \; ,\eqno(8)
$$
where $Q_\beta (m,n)=\sum_{\alpha\in\cA}Q_{\alpha\beta}
(m,n)$ is the return time of $x$ in $I(n)$.

\vskip .3 truecm
\noindent
{\bf 3.2} 
The operators $S (m,n)$ preserve regularity and 
commute with derivation. They also satisfy
$$
\int_{\sqcup_\cA I_\alpha (m)}\varphi 
= \int_{\sqcup_\cA I_\beta (n)} S (m,n)\varphi\; .\eqno(9)
$$
If the restriction of $\varphi$ to each of the intervals 
$I_\alpha (m)$ is a polynomial of degree $\le\mu$ then the restriction 
of $S (m,n)\varphi$ to each of the intervals 
$I_\beta (n)$ is also a polynomial of degree $\le\mu$.

We will denote $\Gamma^{(m)}$ the space of functions 
$\varphi$ which are constant on each of the intervals 
$I_\alpha (m)$. The characteristic functions of the intervals 
$I_\alpha (m)$ form a basis of $\Gamma^{(m)}$. The operator
$S (m,n)$ maps $\Gamma^{(m)}$ into $\Gamma^{(n)}$.
In the bases we have chosen of $\Gamma^{(m)}$ and 
$\Gamma^{(n)}$ the matrix of $S(m,n)|_{\Gamma^{(m)}}$
is ${}^tQ (m,n)$.
We denote $\Gamma^{(m)}_*$ the hyperplane of 
$\Gamma^{(m)}$  whose elements are the functions
whose integral on the disjoint union $\sqcup I_\alpha (m)$ vanishes.
It is sent by $S(m,n)$ into $\Gamma^{(n)}_*$.

\vskip .3 truecm
\noindent
{\bf 3.3}
Write $\hbox{BV}_*\, (m)$ for the space of functions 
$\varphi$ on $\sqcup_\cA I_\alpha (m)$ which are of 
bounded variation on each $I_\alpha (m)$ and of mean value $0$;
write $\hbox{BV}^1\, (m)$ for the space of functions 
$\Phi$ on $\sqcup_\cA I_\alpha (m)$ which are lipschitzian on each 
$I_\alpha (m)$ and whose derivative belongs to $\hbox{BV}_*\, (m)$.

We denote $E^{(m)}_s$
the space of functions $\Phi\in\hbox{BV}^1\, (m)$
such that
$$
\limsup_{n\rightarrow\infty}{\log\Vert S(m,n)\Phi\Vert_\infty
\over\log \Vert Q(m,n)\Vert}<0\; .\eqno(10)
$$
We obviously have $S(m,n)(E^{(m)}_s)\subset E^{(n)}_s$.
We will denote $\Gamma^{(m)}_s$ 
the intersection of $E^{(m)}_s$ with $\Gamma^{(m)}$.

Since $S(m,n)(\Gamma^{(m)}_s )\subset \Gamma^{(n)}_s$
we can consider the induced operator
$$
S(m,n)_\flat \, :\, \Gamma^{(m)}_*/\Gamma^{(m)}_s\rightarrow
\Gamma^{(n)}_*/\Gamma^{(n)}_s\; .\eqno(11)
$$

\vskip .5 truecm
\noindent
{\bf 4. The Diophantine condition}

\vskip .3 truecm\noindent
Here we introduce the notion of i.e.m.\ of ``Roth type''
which gives a full measure class of i.e.m.\ 's which 
satisfy the assumptions of our theorem.

\vskip .3 truecm\noindent
{\bf 4.1} An i.e.m.\ $T$ is said to be of {\it ``Roth type''}
if its accelerated Zorich continued fraction verifies the following 
conditions: 
\item{(a)} For any $\varepsilon >0$, we have 
$\Vert Z(n+1)\Vert \le \Vert Q(0,n)\Vert^\varepsilon$
for all large enough $n$; 
\item{(b)} There exists $\theta >0$ such that 
$\Vert S(0,n)\vert_{\Gamma^{(0)}_*}
\Vert \le \Vert S(0,n)\vert_{\Gamma^{(0)}}
\Vert^{1-\theta}=
\Vert Q(0,n)\Vert^{1-\theta}$ for all large enough $n$; 
\item{(c)} For any $\varepsilon >0$, $m<n$, with $n$ large 
enough, we have  
$\Vert (S(m,n)_\flat )^{-1}\Vert\le\Vert Q(0,n)\Vert^\varepsilon$.

\vskip .3 truecm\noindent
It can be shown that for any admissible pair 
$(\pi_0,\pi_1)$, conditions (a), (b), (c) are satisfied by a 
set $D(\pi_0,\pi_1)$ of full measure in $\Delta (\pi_0,\pi_1)$.

\vskip .5 truecm
\noindent
{\bf 5. Sketch of the proof of the theorem}

\vskip .3 truecm\noindent
{\bf 5.1}
When $\hat{T}$ is a minimal homeomorphism of a compact space $X$, we know
from a theorem of
Gottschalk-Hedlund [GH] that a continuous function $\Phi$ on $X$ is 
a $\hat T$--coboundary of some continuous function as soon as 
the Birkhoff sums of $\Phi$ at some point of $X$ are 
bounded. 

Let $T$ be an i.e.m.\ with no orbit segment starting and ending 
with discontinuities; then $T$ is minimal, but not 
continuous. Nevertheless, a Denjoy--like construction 
allows to apply Gottschalk-Hedlund's theorem and conclude that a 
continuous function whose Birkhoff sums at some point are 
bounded is the $T$--coboundary of a bounded function.
Given $\varphi\in\hbox{BV}_*\, (0)$, it is therefore sufficient
to find a primitive $\Phi$ (determined by $d$ constants of integration)
whose Birkhoff sums at $0$ are bounded.

\vskip .3 truecm\noindent
{\bf 5.2} 
Given $N>0$ and a function $\Phi$ on $\sqcup_\cA I_\alpha$, we 
can write the Birkhoff sums of $\Phi$ at $0$ as a finite sum:
$$
\sum_{i=0}^{N-1} \Phi\circ T^i (0)
=\sum_j S(0,n_j)(\Phi )(x_j)\; ,\eqno(12)
$$
where for every $n\ge 0$ we have 
$$
\hbox{card}\, \{j\, ,\, n_j=n\}
\le \Vert Z(n+1)\Vert := \max_\beta
\sum_\alpha Z_{\alpha\beta} (n+1)\; .\eqno(13)
$$
Thus the Birkhoff sums of $\Phi$ at $0$ are bounded 
as soon as 
$$
\sum_{n\ge 0}\Vert Z(n+1)\Vert
\Vert S(0,n)(\Phi )\Vert_{\hbox{L}^\infty}<+\infty\; .\eqno(14)
$$

\vskip .3 truecm\noindent
{\bf 5.3} Let $\varphi\in\hbox{BV}_*\, (0)$; write $\varphi 
=\varphi_0+\chi_0$, where $\chi_0\in\Gamma^{(0)}_*$ and the 
mean value of $\varphi_0$ on every $I_\alpha$ vanishes. 
Write then inductively 
$$
S(n-1,n)(\varphi_{n-1})=\varphi_n+\chi_n\; ,\eqno(15)
$$
where $\chi_n\in\Gamma^{(n)}_*$ and the 
mean value of $\varphi_n$ on every $I_\alpha (n)$ vanishes.
We have, for every $n\ge 0$:
$$
\Vert \varphi_n\Vert_{\hbox{L}^\infty}\le 
\max_\alpha \hbox{Var}_{I_\alpha (n)}\varphi_n
\le \sum_\alpha \hbox{Var}_{I_\alpha }\varphi
:= \hbox{Var}\, \varphi\; ,\eqno(16)
$$
and thus 
$$
\eqalignno{
\Vert\chi_0\Vert &\le \Vert\varphi\Vert_{\hbox{L}^\infty}\; , &(17)\cr
\Vert\chi_n\Vert_{\hbox{L}^\infty}&\le 
\Vert S(n-1,n)(\varphi_{n-1})\Vert_{\hbox{L}^\infty}
\le \Vert Z(n)\Vert\hbox{Var}\, \varphi\; , &(18) \cr
}
$$
for $n>0$. From this and conditions (a), (b) in the definition of 
Roth type, we get easily that there exists $\theta'=\theta' (\theta )>0$
such that 
$$
\Vert S(0,n)(\varphi )\Vert_{\hbox{L}^\infty}\le
\Vert Q(0,n)\Vert^{1-\theta'}\Vert\varphi\Vert_{\hbox{BV}}
\eqno(19)
$$
for all $n\ge 0$ large enough.

\vskip .3 truecm\noindent
{\bf 5.4} Given $\varphi\in
\hbox{BV}_*\, (0)$, we will find a primitive $\Phi\in
\hbox{BV}_*^1\, (0)$ that satisfies (14);
for this, in view of condition (a) in Section 4.1, it is sufficient
to find $\Phi\in E^{(0)}_s$. Thus, we only need to define $\Phi$ 
mod $\Gamma_s^{(0)}$.

For $n\ge 0$ and $\varphi\in  \hbox{BV}_*\, (n)$, let 
$P_0^{(n)}\varphi\in \hbox{BV}_*^1\, (n)$ be the primitive of 
$\varphi$ whose mean value over every interval $I_\alpha (n)$
vanishes; this is not functorial with respect to special Birkhoff 
sums. We define, for $n>0$:
$$
\eqalign{
\Lambda^{(n)} \, &:\, \hbox{BV}_*\, (n-1)\rightarrow
\Gamma_*^{(n)}/\Gamma_s^{(n)}\cr
\Lambda^{(n)} &= P_0^{(n)}S(n-1,n)-S(n-1,n)P_0^{(n-1)}
\hbox{mod}\, \Gamma_s^{(n)}\; .\cr}\eqno(20)
$$
From condition (a) in Section 4.1, it is not difficult to show, 
for every $\varepsilon >0$, $\varphi\in \hbox{BV}_*\, (n-1)$,
$n$ large enough that:
$$
\Vert \Lambda^{(n)}\varphi \Vert_{\hbox{L}^\infty}\le
\Vert Q(0,n)\Vert^{-1+\varepsilon}\Vert\varphi \Vert_{\hbox{L}^\infty}
\; .\eqno(21)
$$
Joining (19), (21) and condition (c) of
Section 4.1, we see that the series 
$$
\Delta P^{(m)}(\varphi ) =\sum_{n>m}\left((S(m,n)_\flat)^{-1}
\circ\Lambda^{(n)}\circ S(m,n-1)\right)(\varphi )\eqno(22)
$$
for $m\ge 0$, $\varphi\in \hbox{BV}_*\, (m)$ converges in 
$\Gamma_*^{(m)}/\Gamma_s^{(m)}$.

\noindent
Setting 
$$
P^{(m)}(\varphi )=P_0^{(m)}(\varphi )+\Delta P^{(m)}(\varphi )
\eqno(23)
$$
we have now functoriality:
$$
S(m,n)(P^{(m)}(\varphi ))=P^{(n)}(S(m,n)\varphi )
\hbox{mod}\, \Gamma_s^{(n)}\; .\eqno(24)
$$
It is now not difficult to check that given $\varphi\in 
\hbox{BV}_*\, (0)$, any primitive $\Phi\in \hbox{BV}_*^1\, (0)$
whose class $\hbox{mod}\, \Gamma_s^{(0)}$ belongs to 
$P^{(0)}(\varphi )$ also belongs to 
$E_s^{(0)}$, and thus is the coboundary of a bounded function.

\vskip .3 truecm
\noindent
{\bf References}

\vskip .3 truecm
\item{[Fo]} G. Forni ``Solutions of the cohomological equation
for area-preserving flows on compact surfaces of higher genus''
{\it Annals of Mathematics} {\bf 146} (1997)
295-344.
\item{[GH]} W.H. Gottschalk, G.A. Hedlund ``Topological 
Dynamics'' A.M.S. Coll. Publ. {\bf 36} (1955)
\item{[Ke1]} M. Keane ``Interval exchange transformations''
{\it Math. Z.} {\bf 141} (1975) 25--31
\item{[Ke2]} M. Keane ``Non--ergodic interval exchange transformations''
{\it Isr. J. Math.} {\bf 26} (1977) 188--196
\item{[KN]} H. B. Keynes and D. Newton ``A ``Minimal'',
Non--Uniquely Ergodic Interval Exchange Transformation''
{\it Math. Z.} {\bf 148} (1976) 101--105
\item{[Ma]} H. Masur ``Interval exchange transformations and measured
foliations'' Annals of Mathematics 115 (1982) 169--200
\item{[Ra]} G. Rauzy ``\'Echanges d'intervalles et transformations
induites'' {\it Acta Arit.} (1979) 315--328
\item{[V1]} W. Veech ``Interval exchange transformations''
Journal d'Analyse Ma\-th\'e\-ma\-ti\-que 33 (1978) 222-272
\item{[V2]} W. Veech ``Gauss measures for transformations on the space of
interval exchange maps'' Ann. of Math. 115 (1982) 201--242
\item{[Z1]} A. Zorich ``Finite Gauss measure on the space of interval
exchange transformations. Lyapunov exponents''
Annales de l'Institut Fourier Tome
46 fasc. 2 (1996) 325-370
\item{[Z2]} A. Zorich ``Deviation for interval exchange
transformations'' {\it Ergod. Th. Dyn. Sys.}{\bf 17} (1997), 1477--1499

\bye